# Trigonometry and Analytic Tools in Olympiad Geometry Problems, Part I

Orestis Lignos (orelig2006@gmail.com)

December 2023

## Contents

# 1 Abstract

This is the first part of a series of papers aiming to show how trigonometry and analytic tools can help into tackling demanding Olympiad geometry problems. We present several novel techniques for tackling hard problems from various international and national mathematical competitions, and develop the appropriate theoretical background.

Several insightful examples and practice problems are included for the reader to sharpen their problem solving skills. The number of asterisks next to each problem indicate its (subjective) difficulty level, ranging from **(*)** (early IMO problem) to **(***)** (hard IMO problem), or ocassionaly **(****)** (too hard for an IMO problem). Lastly, at the end of the paper we devote a section on providing some historical information on mathematicians whose names appeared in the paper.

I would like to point out that no diagrams are provided despite this bring a geometry handout. This was done on purpose; trigonometry solutions do not focus on exploring the diagram. Instead, they focus on eliciting crucial information from it, as we will see below. The diagrams for the majority of the problems presented can even be drawn by hand.

I am thankful to Konstantinos Konstantinidis for bringing some problems into my attention, and I dedicate this series of papers to my parents, for being my beloved family, and to the Hellenic Mathematical Society and the Greek Mathematical Olympiads Community, for being my second family for so many years.

## 2 The main weapons, Part 1

Since we are going to discuss how to solve geometry problems with trigonometry, we should firstly state the following two simple facts (throughout the whole paper we denote the sides $AB, BC, CA$ of a triangle $ABC$ by $c, a, b$ respectively):

> **Theorem 1: Law of Sines**
>
> Let $ABC$ be a triangle, and $R$ be the radius of its circumcircle. Then,
> $$\frac{a}{\sin \angle A} = \frac{b}{\sin \angle B} = \frac{c}{\sin \angle C} = 2R$$

*Proof:* Omitted ∎

Moving on, here is the next one:

> **Theorem 2: Law of Cosines**
>
> Let $ABC$ be a triangle, and $R$ be the radius of its circumcircle. Then,
> $$a^2 = b^2 + c^2 - 2bc \cos \angle A,$$
> and we can be happy since similar relations hold for $\cos \angle B$ and $\cos \angle C$, too.

*Proof:* Omitted ∎

Some fundamental and useful trigonometric identities follow, in order to enrich our toolbox.

> **Theorem 3: Trigonometric identities**
>
> We have the following list of identities:
>
> - $\sin(2x) = 2 \sin x \cos x$
> - $\cos(2x) = 2 \cos^2 x - 1 = 1 - 2 \sin^2 x$
> - $\sin(x + y) = \sin x \cos y + \cos x \sin y$
> - $\sin(x - y) = \sin x \cos y - \cos x \sin y$
> - $\cos(x + y) = \cos x \cos y - \sin x \sin y$
> - $\cos(x - y) = \cos x \cos y + \sin x \sin y$

*Proof:* Omitted ∎

**Exercise 1:** Prove the identity $\cos(x+y) = \cos x \cos y - \sin x \sin y$ (and proceed similarly for each one of the last four identities) using Euler's formula $e^{i\theta} = \cos\theta + i\sin\theta$. Can you give the respective geometric interpretation?

What follows is an important tool for our toolbox:

> **Lemma 1: Ratio Lemma**
>
> Let $ABC$ be a triangle, and $M$ be a point on $BC$. Then,
> $$\frac{BM}{MC} = \frac{AB}{AC} \cdot \frac{\sin \angle BAM}{\sin \angle MAC}$$

*Proof:* Apply the Law of Sines in triangles $ABM$ and $ACM$, to get:
$$\frac{BM}{\sin \angle BAM} = \frac{AB}{\sin \angle AMB}$$
and
$$\frac{CM}{\sin \angle CAM} = \frac{AC}{\sin \angle AMC}$$

Now note that $\sin \angle AMB = \sin \angle AMC$ since they sum to $\pi$, divide these two relations and we are done ∎

Proving the Butterfly theorem with the Ratio Lemma is a particularly instructive ans straightforward exercise, which we leave to the reader:

**Exercise 2 (Butterfly theorem):** Let $M$ be the midpoint of a chord $PQ$ of a circle, through which two other chords $AB$ and $CD$ are drawn. $AD$ and $BC$ intersect chord $PQ$ at $X$ and $Y$ respectively. Prove that point $M$ is the midpoint of $XY$ in two ways, using

  (i) the Ratio lemma and

  (ii) [Haruki's Lemma](#).

Here is another Lemma, which we present here since its proof uses the Law of Cosines (and it is very useful in length bashing proofs):

> **Lemma 2: Perpendicularity lemma**
>
> Let $A, B, C, D$ be points in the plain. Then,
> $$AC \perp BD \iff AB^2 - AD^2 = BC^2 - CD^2$$

*Proof:* The $\implies$ direction is merely some easy applications of the Pythagorean Theorem. For the $\impliedby$ direction, let $\angle AXB = \angle CXD = \omega$ and $X \equiv AC \cap BD$ and so by the Law of Cosines
$$AB^2 = AX^2 + XB^2 - 2AX \cdot XB \cos \angle \omega$$
and
$$AD^2 = AX^2 + XD^2 + 2AX \cdot XD \cos \angle \omega,$$
hence
$$AB^2 - AD^2 = XB^2 - XD^2 - 2AX \cdot BD \cos \angle \omega$$
and similarly
$$CB^2 - CD^2 = XB^2 - XD^2 + 2CX \cdot BD \cos \angle \omega$$
Thus,
$$AB^2 - AD^2 = BC^2 - CD^2 \iff 2AC \cdot BD \cos \angle \omega = 0,$$
and so $\cos \angle \omega = 0$ i.e. $\omega = 90°$, as desired ∎

Here is an easy exercise illustrating the importance and efficiency of this Lemma:

**Exercise 3:** Let $I$ be the center of the incircle and $D$ be the antipode of $A$ in triangle $ABC$. Points $E$ and $F$ belong on lines $BA$ and $CA$ such that $BE = CF = \dfrac{AB + BC + CA}{2}$. Prove that $EF$ is perpendicular to $OI$.

Here are some examples which come from various competitions that nicely illustrate the importance of the above results.

> **Example 1: Balkan MO 2022/1**
>
> Let $ABC$ be an acute triangle such that $CA \neq CB$ with circumcircle $\omega$ and circumcentre $O$. Let $t_A$ and $t_B$ be the tangents to $\omega$ at $A$ and $B$ respectively, which meet at $X$. Let $Y$ be the foot of the perpendicular from $O$ onto the line segment $CX$. The line through $C$ parallel to line $AB$ meets $t_A$ at $Z$. Prove that the line $YZ$ passes through the midpoint of the line segment $AC$.

*Proof:* This problem is a typical triangle geometry, and we recognize it's susceptible to bash. Let $ZY$ intersect $AC$ at $M$.

From the Ratio Lemma we know that
$$\frac{AM}{MC} = \frac{AZ}{ZC} \cdot \frac{\sin \angle AZM}{\sin \angle MZC}$$
and
$$\frac{CY}{YX} = \frac{CZ}{ZX} \cdot \frac{\sin \angle CZY}{\sin \angle YZX}$$

and so by multiplying these two we get
$$\frac{AM}{MC} \cdot \frac{CY}{YX} = \frac{AZ}{ZX}$$
(note that this relation also follows from a straightforward application of Menelaus's theorem, which will be developed in a later section), hence in order to conclude we need to prove that
$$\frac{CY}{YX} = \frac{AZ}{ZX}$$

Now, let $CX$ meet the circumcircle of $ABC$ at $T$ and $AB$ at $D$. By Thales, we know that
$$\frac{AZ}{ZX} = \frac{CD}{CX}$$
and so we are left to prove
$$\frac{CD}{CX} = \frac{CY}{YX}$$
Now, subtracting the numerator from the denominator we are left with
$$\frac{CD}{DX} = \frac{CY}{TX}$$

And now, time for the Ratio Lemma once again! We have that
$$\frac{CD}{DX} = \frac{CA}{AX} \cdot \frac{\sin \angle A}{\sin \angle BAX} = \frac{b}{AX} \cdot \frac{a}{c} = \frac{ab}{c \cdot AX}$$
and
$$\frac{CY}{TX} = \frac{CT}{2TX} = \frac{1}{2} \cdot \frac{CA}{XA} \cdot \frac{\sin \angle CAT}{\sin \angle TAX} = \frac{1}{2} \cdot \frac{b}{XA} \cdot \frac{\sin \angle CAT}{\sin \angle ACT} = \frac{b}{XA} \cdot \frac{CT}{2AT}$$

So, in order to conclude we just need to show that
$$\frac{CT}{2AT} = \frac{a}{c},$$
or equivalently that triangles $CAT$ and $CBN$ are similar, where $N$ is the midpoint of $AB$. However, this is a well-known property of the symmedian (which we will prove in the next section for those that have not seen it yet), and so we are done ∎

For comparison reasons, we present the synthetic solution for this problem, too

*Proof 2 (synthetic):* Note that $\angle OEX = \angle OAX = \angle OBX = 90°$, and so $O, Y, A, X, B$ are concyclic. Thus,
$$\angle AYX = \angle ABX = \angle BAX = \angle CZA,$$
and so $ZCYA$ is cyclic. Therefore
$$\angle CYZ = \angle CAZ = \angle B$$
Now, define $M$ to be the midpoint of $CA$. Thus,
$$\angle CMO = \angle CYO = 90°,$$
which implies that $CMYO$ is cyclic, and so
$$\angle CYM = \angle COM = \angle B$$
To sum up,
$$\angle CYM = \angle B = \angle CYZ,$$
so $Z, M, Y$ are collinear, as desired ∎

We move on to our next example:

> **Example 2: Adapted from the Internet**
>
> Let $ABC$ be a triangle, and let its incircle with center $I$ meet lines $AB, AC$ at $F, E$ respectively. If $M$ is the midpoint of $BC$ and $AM$ meets $FE$ at $N$, then prove that $NI$ is perpendicular to $BC$.

*Proof:* Suppose the incircle touches $BC$ at point $D$, and define $N$ as the point where lines $DI$ and $FE$ meet. Then we just need to prove that $A, N, M$ are collinear. Note that by the Ratio Lemma
$$\frac{FN}{NE} = \frac{FD}{DE} \cdot \frac{\sin \angle FDN}{\sin \angle NDE} =$$
$$= \frac{\sin \angle FED}{\sin \angle EFD} \cdot \frac{\sin \angle ABI}{\sin \angle ACI} = \frac{\cos \angle B/2}{\cos \angle C/2} \cdot \frac{\sin \angle B/2}{\sin \angle C/2} =$$
$$= \frac{\sin \angle B}{\sin \angle C} = \frac{AC}{AB}$$
Moreover,
$$\frac{FN}{NE} = \frac{FA}{AE} \cdot \frac{\sin \angle FAN}{\sin \angle NAE}$$
and since $FA = AE$, we finally obtain
$$\frac{\sin \angle FAN}{\sin \angle NAE} = \frac{AC}{AB}$$

Now let, $AN$ meet $BC$ at $M'$. Then, by the Ratio Lemma

$$\frac{BM'}{M'C} = \frac{AB}{AC} \cdot \frac{\sin \angle FAN}{\sin \angle NAE} = 1,$$

and so $M'$ is the midpoint of $BC$, as desired ∎

Our next example comes from Bulgaria, and it has an interesting story behind it.

> **Example 3: Bulgaria 2021/2**
>
> A point $T$ is given on the altitude through point $C$ in the acute triangle $ABC$ with circumcenter $O$, such that $\angle TBA = \angle ACB$. If the line $CO$ intersects side $AB$ at point $K$, prove that the perpendicular bisector of $AB$, the altitude through $A$ and the segment $KT$ are concurrent.

*Proof:* Let $M$ be the midpoint of $AB$, $D$ be the foot of the altitude from $C$, and suppose that $OM$ meets the altitude $AS$ at point $P$. We need to show that $K, P, T$ are collinear. Since $OM$ and $DC$ are parallel (they are both perpendicular to $AB$), this is equivalent to proving that

$$\frac{MP}{PO} = \frac{DT}{TC}$$

Note that (we assume that $\angle B > \angle C$)

$$\angle MAP = 90° - \angle B \text{ and } \angle PAO = \angle BAO - \angle BAP = \angle B - \angle C$$

and so by the Ratio Lemma,

$$\frac{MP}{PO} = \frac{MA}{AO} \cdot \frac{\sin \angle MAP}{\sin \angle PAO} = \sin \angle AOM \cdot \frac{\cos \angle B}{\sin(\angle B - \angle C)} = \frac{\sin \angle C \cos \angle B}{\sin(\angle B - \angle C)}$$

Moreover, again by the Ratio Lemma,

$$\frac{DT}{TC} = \frac{DB}{BC} \cdot \frac{\sin \angle DBT}{\sin \angle TBC} = \cos \angle B \cdot \frac{\sin \angle C}{\sin(\angle B - \angle C)}$$

So the two ratios are equal, and we are done ∎

This problem still seeks a synthetic solution in the respective Art of Problem Solving thread. Using our tools though, it becomes routine.

We move on to a recent IMO problem.

> **Example 4: IMO 2022/4**
>
> Let $ABCDE$ be a convex pentagon such that $BC = DE$. Assume that there is a point $T$ inside $ABCDE$ with
>
> $$TB = TD, TC = TE \text{ and } \angle ABT = \angle TEA$$
>
> Let line $AB$ intersect lines $CD$ and $CT$ at points $P$ and $Q$, respectively. Assume that the points $P, B, A, Q$ occur on their line in that order.
>
> Let line $AE$ intersect $CD$ and $DT$ at points $R$ and $S$, respectively. Assume that the points $R, E, A, S$ occur on their line in that order.
>
> Prove that the points $P, S, Q, R$ lie on a circle.

*Proof:* The triangles $BTC$ and $ETD$ are equal, since they have three equal sides. Therefore, $\angle BTC = \angle ETD$. Hence,

$$\angle AQT = \angle BTC - \angle ABT = \angle ETD - \angle AET = \angle AST$$

Therefore, by the Law of Sines in triangles $AST, AQT$ we get

$$\frac{ST}{\sin \angle SAT} = \frac{AT}{\sin \angle AST} = \frac{AT}{\sin \angle AQT} = \frac{QT}{\sin \angle QAT}$$

thus

$$\frac{ST}{\sin \angle SAT} = \frac{QT}{\sin \angle QAT}$$

which implies that $\dfrac{ST}{QT} = \dfrac{\sin \angle SAT}{\sin \angle QAT}$ **(1)**

Moreover, by the Law of Sines in triangles $ABT, AET$ we get

$$\frac{BT}{\sin \angle BAT} = \frac{AT}{\sin \angle ABT} = \frac{AT}{\sin \angle AET} = \frac{ET}{\sin \angle EAT}$$

thus

$$\frac{BT}{\sin \angle BAT} = \frac{ET}{\sin \angle EAT}$$

which implies that $\dfrac{TD}{TC} = \dfrac{BT}{ET} = \dfrac{\sin \angle BAT}{\sin \angle EAT} = \dfrac{\sin \angle QAT}{\sin \angle SAT}$, and so $\dfrac{TC}{TD} = \dfrac{\sin \angle SAT}{\sin \angle QAT}$ **(2)**

By (1) and (2) we get that $\dfrac{ST}{QT} = \dfrac{TC}{TD}$, hence the quadrilateral $QSCD$ is cyclic. Therefore, $\angle QSD = \angle QCD$. Hence,

$$\angle QSR = \angle QSD - \angle AST = \angle QCD - \angle AQT = \angle QPR,$$

which gives us that ∠$QSR$ = ∠$QPR$, hence we conclude that $QSPR$ is cyclic, as desired ∎

Our next example is a very useful lemma.

> **Example 5: Steiner's relation**
>
> Let $ABC$ be a triangle and $D, E$ be points on $BC$ such that $AD, AE$ are isogonal with regard to angle $A$. Then,
> $$\frac{BD}{DC} \cdot \frac{BE}{EC} = \frac{AB^2}{AC^2}$$

*Proof:* Apply the Ratio Lemma twice, then multiply the two resulting relations! ∎

As a corollary of this relation, we obtain the following Lemma (also refer to Symmedian lemma II in the next section):

> **Lemma 3: Symmedian Lemma I**
>
> If $ABC$ is a triangle and $T$ a point on $BC$ such that $AT$ is the $A-$ symmedian, then
> $$\frac{BT}{TC} = \frac{AB^2}{AC^2}$$

*Proof:* By definition, $AT$ is isogonal with the median $AM$, and now we may apply Steiner's relation ∎

Our next problem comes from USA TSTST, and is a very instructive one.

> **Example 6: USA TSTST 2019**
>
> Let $ABC$ be an acute triangle with orthocenter $H$ and circumcircle $\Gamma$. A line through $H$ intersects segments $AB$ and $AC$ at $E$ and $F$, respectively. Let $K$ be the circumcenter of $\triangle AEF$, and suppose line $AK$ intersects $\Gamma$ again at a point $D$. Prove that line $HK$ and the line through $D$ perpendicular to $\overline{BC}$ meet on $\Gamma$.

*Proof (due to AoPS member Eyed):* Let $Q$ be on $\Gamma$ such that $DQ \perp BC, P = DQ \cap BC, X = DQ \cap EF, R = AK \cap EF$, and $H'$ the reflection of $H$ over $BC$ (which lies on $\Gamma$). Let ∠$AFE$ = $\beta$, ∠$AEF$ = $\gamma$. We claim that $AQ$ is parallel to $EF$. This is because
$$∠AQX = ∠C + ∠BCD = ∠C + 90 - \beta$$

while
$$\angle EXP = 360 - (90 + \angle B + 180 - \gamma) = 90 - \angle B + \gamma = 90 + \angle C - \beta = \angle AQX$$

Thus, those two are parallel. Since $AH$ and $QX$ are perpendicular to $BC$, we have that $AQHX$ is a parallelogram, so $AQ = HX$. In addition, since $AH'$ is parallel to $QD$, this means that $AQ = H'D = HX$, which implies $X$ is the reflection of $D$ over $BC$ and thus $BHXC$ is cyclic.

Now, let $K'$ be the intersection of $HQ$ and $AD$. Observe that
$$\frac{AK}{KR} = \frac{EK}{KR} = \frac{\sin \angle KRE}{\sin \angle KER} = \frac{\sin(90 + \beta - \gamma)}{\sin(90 - A)} = \frac{\cos(\beta - \gamma)}{\cos A}$$

Now, we will calculate $\dfrac{AK'}{K'R}$. Since $HR$ is parallel to $AQ$ and $AH$ is parallel to $QX$, we have:
$$\frac{AK'}{K'R} = \frac{K'Q}{K'H} = \frac{QD}{AH} = \frac{QD}{QX}$$

We can now directly calculate $QD$ and $QX$. By the Law of Sines on triangles $QDC, QXC$, we have
$$\frac{QD}{DC} = \frac{\sin \angle DCQ}{\sin \angle DQC} = \frac{\sin(\angle BDP + \angle BAD)}{\sin(90 - \gamma)} = \frac{\sin(\gamma + 90 - \beta)}{\cos(\gamma)} = \frac{\cos(\gamma - \beta)}{\cos(\gamma)}$$

and
$$\frac{QX}{DC} = \frac{QX}{XC} = \frac{\sin \angle XCQ}{\sin \angle XQC} = \frac{\sin \angle DCQ - 2(90 - \beta)}{\sin(90 - \gamma)} = \frac{\sin(\gamma + 90 - \beta - 180 + 2\beta)}{\cos(\gamma)} = -\frac{\cos(\beta + \gamma)}{\cos \gamma}$$

Thus,
$$\frac{QD}{QX} = \frac{\frac{QD}{DC}}{\frac{QX}{DC}} = \frac{\frac{\cos(\gamma - \beta)}{\cos \gamma}}{-\frac{\cos(\beta + \gamma)}{\cos \gamma}} = -\frac{\cos(\beta - \gamma)}{\cos(180 - A)} = \frac{\cos(\beta - \gamma)}{\cos A} = \frac{AK'}{K'R}$$

Since $\dfrac{AK'}{K'R} = \dfrac{AK}{KR}$, and they both lie on $AR$, we conclude that $K' = K$, which proves the problem ∎

We end this section with a last example, coming from the European's Girls Mathematical Olympiad (EGMO).

> **Example 7: EGMO 2021/3**
>
> Let $ABC$ be a triangle with an obtuse angle at $A$. Let $E$ and $F$ be the intersections of the external bisector of angle $A$ with the altitudes of $ABC$ through $B$ and $C$ respectively. Let $M$ and $N$ be the points on the segments $EC$ and $FB$ respectively such that $\angle EMA = \angle BCA$ and $\angle ANF = \angle ABC$. Prove that the points $E, F, N, M$ lie on a circle.

*Proof:* Let $H$ be the orthocenter of triangle $ABC$. We will, in fact, prove that points $E, F, N, M, H$ are all concyclic, which will yield the desired result. To this end, it suffices to prove that $E, N, F, H$ are concyclic since then by similar arguments we would have that $E, M, F, H$ are concylic too. We will, therefore, prove that $BE \cdot BH = BN \cdot BF$

By the Law of Sines in triangle $ABN$,

$$\frac{BN}{\sin \angle BAN} = \frac{BA}{\sin \angle ANB} = \frac{c}{\sin \angle B}$$

and since $\angle BAN = \angle B - \angle ABN = \angle FBC$, we obtain that

$$BN = \frac{c \sin \angle FBC}{\sin \angle B} \quad \textcolor{red}{(1)}$$

By the Law of Sines in triangle $BFC$,

$$\frac{FC}{\sin \angle FBC} = \frac{FB}{\sin \angle HCB}$$

and so

$$BF = \frac{FC \sin \angle HCB}{\sin \angle FBC} \quad \textcolor{red}{(2)}$$

Thus, multiplying (1) and (2) we obtain

$$BN \cdot BF = \frac{c \cdot FC \sin \angle HCB}{\sin \angle B}$$

and so it suffices to prove that

$$\frac{c \cdot FC \sin \angle HCB}{\sin \angle B} = BE \cdot BH$$

which rearranges to

$$\frac{c}{\sin \angle B} \cdot \frac{\sin \angle HCB}{BH} = \frac{BE}{CF} \quad \textcolor{red}{(3)}$$

Now, we do the following two computations to end the problem:

- By the Law of Sines in triangle $HBC$,

$$\frac{BH}{\sin \angle HCB} = \frac{BC}{\sin \angle BHC} = \frac{BC}{\sin \angle A} = \frac{b}{\sin \angle B}$$

and so the left hand side of (3) is equal to $\frac{c}{b}$ and

- Note that $\angle EAB = \angle FAC = 90° - \angle A/2$ and by an easy angle-chase $\angle AEB = \angle AFC = 180° - \angle A/2$. Hence, by the Law of Sines in triangles $ABE, ACF$,

$$\frac{BE}{\sin \angle EAB} = \frac{AB}{\sin \angle AEB} \text{ and } \frac{CF}{\sin \angle FAC} = \frac{AC}{\sin \angle AFC}$$

which when divided give us $\dfrac{BE}{CF} = \dfrac{c}{b}$ and so the right hand side of (3) is also equal to $\dfrac{c}{b}$, as desired.

Hence, the problem is solved ∎

# 3  Practice Problems

**Problem 1: Balkan MO 2020/1 (*)**

Let $ABC$ be an acute triangle with $AB = AC$, let $D$ be the midpoint of the side $AC$, and let $\gamma$ be the circumcircle of the triangle $ABD$. The tangent of $\gamma$ at $A$ crosses the line $BC$ at $E$. Let $O$ be the circumcenter of the triangle $ABE$. Prove that midpoint of the segment $AO$ lies on $\gamma$.

**Problem 2: IMO Shortlist 2016 G2 (*)**

Let $ABC$ be a triangle with circumcircle $\Gamma$ and incenter $I$ and let $M$ be the midpoint of $\overline{BC}$. The points $D$, $E$, $F$ are selected on sides $\overline{BC}$, $\overline{CA}$, $\overline{AB}$ such that $\overline{ID} \perp \overline{BC}$, $\overline{IE} \perp \overline{AI}$, and $\overline{IF} \perp \overline{AI}$. Suppose that the circumcircle of $\triangle AEF$ intersects $\Gamma$ at a point $X$ other than $A$. Prove that lines $XD$ and $AM$ meet on $\Gamma$.

**Problem 3: IMO Shortlist 2020 G1 (*)**

Let $ABC$ be an isosceles triangle with $BC = CA$, and let $D$ be a point inside side $AB$ such that $AD < DB$. Let $P$ and $Q$ be two points inside sides $BC$ and $CA$, respectively, such that $\angle DPB = \angle DQA = 90°$. Let the perpendicular bisector of $PQ$ meet line segment $CQ$ at $E$, and let the circumcircles of triangles $ABC$ and $CPQ$ meet again at point $F$, different from $C$. Suppose that $P$, $E$, $F$ are collinear. Prove that $\angle ACB = 90°$.

**Problem 4: Balkan MO 2013/1 (**)**

In a triangle $ABC$, the excircle $\omega_a$ opposite $A$ touches $AB$ at $P$ and $AC$ at $Q$, while the excircle $\omega_b$ opposite $B$ touches $BA$ at $M$ and $BC$ at $N$. Let $K$ be the projection of $C$ onto $MN$ and let $L$ be the projection of $C$ onto $PQ$. Show that the quadrilateral $MKLP$ is cyclic.

**Problem 5: IMO 2009/4 (**)**

Let $ABC$ be a triangle with $AB = AC$. The angle bisectors of $\angle CAB$ and $\angle ABC$ meet the sides $BC$ and $CA$ at $D$ and $E$, respectively. Let $K$ be the incentre of triangle $ADC$. Suppose that $\angle BEK = 45°$. Find all possible values of $\angle CAB$.

**Problem 6: USAMO 2019/2 (\*\*)**

Let $ABCD$ be a cyclic quadrilateral satisfying $AD^2 + BC^2 = AB^2$. The diagonals of $ABCD$ intersect at $E$. Let $P$ be a point on side $\overline{AB}$ satisfying $\angle APD = \angle BPC$. Show that line $PE$ bisects $\overline{CD}$.

**Problem 7: China TST 2005 (\*\*)**

Let $ABC$ be a triangle, and let $(\gamma)$ be a circle tangent to $AB, AC$ at $P, Q$ and to the circumcircle of triangle $ABC$ internally at $T$. Suppose that $AT, PQ$ meet at $S$. Prove that $\angle ABS = \angle ACS$.

**Problem 8: USAMO 2014/5 (\*\*)**

Let $ABC$ be a triangle with orthocenter $H$ and let $P$ be the second intersection of the circumcircle of triangle $AHC$ with the internal bisector of the angle $\angle BAC$. Let $X$ be the circumcenter of triangle $APB$ and $Y$ the orthocenter of triangle $APC$. Prove that the length of segment $XY$ is equal to the circumradius of triangle $ABC$.

**Problem 9: Own and Taes Padhihary, OTJMO 2021 (\*\*)**

In triangle $ABC$ with circumcircle $\Gamma$, let $\ell_1$, $\ell_2$, and $\ell_3$ be the tangents to $\Gamma$ at points $A$, $B$, and $C$, respectively. Choose a variable point $P$ on side $\overline{BC}$. Let the lines parallel to $\ell_2$ and $\ell_3$, passing through $P$, meet $\ell_1$ at points $C_1$ and $B_1$, respectively. Let the circumcircles of $\triangle PBB_1$ and $\triangle PCC_1$ meet each other again at a point $Q \neq P$. Let lines $\ell_1$ and $BC$ meet at a point $R$, and let lines $\ell_2$ and $\ell_3$ meet at a point $X$. Prove that, as $P$ varies on side $\overline{BC}$, lines $PQ$ and $RX$ meet at a fixed point.

**Problem 10: CJMO 2023 (\*\*)**

Let $ABC$ be an isosceles triangle with $AB = AC$, and denote by $H$ the orthocenter of $ABC$. Define $N$ as the midpoint of arc $AHC$ (on the circumcircle of $\triangle AHC$), and let $D$ be a point on $\overline{AB}$ such that $\overline{DN}$ is parallel to $\overline{AC}$. If $E$ is the reflection of $D$ across $\overline{AH}$, then show that $\angle DNB$ and $\angle ENC$ are either equal or supplementary.

**Problem 11: GGG4 (***)**

Let $ABC$ be a scalene triangle with incenter $I$ and circumcircle $\Omega$; $L$ is the midpoint of arc $BAC$ and $A'$ is diametrically opposite $A$ on $\Omega$. $D$ is the foot of the perpendicular from $I$ to $\overline{BC}$. $\overleftrightarrow{LI}$ meets $\overleftrightarrow{BC}$ and $\Omega$ at $X$ and $Y$, respectively, and $\overleftrightarrow{LD}$ meets $\Omega$ again at $Z$. $\overleftrightarrow{XZ}$ meets $\overleftrightarrow{A'I}$ at $T$. $\overleftrightarrow{DI}$ meets the circle with diameter $\overline{AI}$ again at $P$. Show that the second intersection between $\overleftrightarrow{PT}$ and the circle with diameter $\overline{AI}$ lies on $\overline{AY}$.

**Problem 12: Vietnam TST 2023 (****)**

Let $ABC$ be an acute, non-isosceles triangle with circumcircle $(O)$. $BE, CF$ are the heights of $\triangle ABC$, and $BE, CF$ intersect at $H$. Let $M$ be the midpoint of $AH$, and $K$ be the point on $EF$ such that $HK \perp EF$. A line not going through $A$ and parallel to $BC$ intersects the minor arc $AB$ and $AC$ of $(O)$ at $P, Q$, respectively.

Show that the tangent line of $(CQE)$ at $E$, the tangent line of $(BPF)$ at $F$, and $MK$ concur.

# 4 The main weapons, Part 2

Here is our next important Lemma:

**Lemma 4: The $\alpha + \gamma = \beta + \delta$ lemma**

If $\alpha, \beta, \gamma, \delta$ are angles such that

$$\alpha + \gamma = \beta + \delta < \pi \text{ and } \frac{\sin \alpha}{\sin \gamma} = \frac{\sin \beta}{\sin \delta}$$

then $\alpha = \beta$ and $\gamma = \delta$

*First proof:* Let $\alpha + \gamma = \beta + \delta = t < \pi$ and consider the function $f(x) = \dfrac{\sin x}{\sin(t-x)}$. This function has derevative

$$f'(x) = \frac{\cos x \sin(t-x) + \sin x \cos(t-x)}{\sin(t-x)^2} = \frac{\sin t}{\sin(t-x)^2},$$

which is positive as $0 < t < \pi$, and so $f$ is strictly increasing, hence injective, and so $f(\alpha) = f(\beta)$ implies $\alpha = \beta$, which also gives $\gamma = \delta$ as desired ∎

*Second proof (due to Demetres Christophides:)* Note that $\sin \alpha \sin \delta = \sin \beta \sin \gamma$, and so
$$\cos(\alpha - \delta) - \cos(\alpha + \delta) = \cos(\beta - \gamma) - \cos(\beta + \gamma),$$
and since $\alpha - \delta = \beta - \gamma$, we infer that $\cos(\alpha + \delta) = \cos(\beta + \gamma)$, which implies that $\alpha + \delta = \beta + \gamma$, and from here we may easily conclude ∎

**Exercise 4:** Try to find a geometric proof of this Lemma!

Now that we have this Lemma, we will come back to what we left unproved in the solution of Example 1.

**Lemma 5: Symmedian lemma II**

Let $ABC$ be a triangle, and suppose the tangents at the circumcircle of this triangle at points $B, C$ meet at point $X$. Let $M$ be the midpoint of $BC$. Then, $AX, AM$ are *isogonal* with regard to angle $A$.

*Proof:* This lemma is fairly well known, but here we aim to provide another proof, using Lemma 4 from above. As we see, this proof requires little or no ingenuity.

In order to prove the statement, we need to prove that $\angle BAX = \angle CAM$ and $\angle CAX = \angle BAM$. Note that we already have the first condition of our Lemma:

$$\angle BAX + \angle CAX = \angle BAM + \angle CAM = \angle A < \pi,$$

and so we just need to verify the second one. Apply the Law of Sines in triangles $ABX$ and $ACX$ to obtain

$$\frac{BX}{\sin \angle BAX} = \frac{AX}{\sin \angle ABX} \text{ and } \frac{CX}{\sin \angle CAX} = \frac{AX}{\sin \angle ACX}$$

Dividing these two we obtain that

$$\frac{\sin \angle CAX}{\sin \angle BAX} = \frac{\sin \angle ACX}{\sin \angle ABX} = \frac{\sin \angle B}{\sin \angle C} = \frac{AC}{AB}$$

Moreover, by the Ratio Lemma

$$\frac{BM}{CM} = \frac{AB}{AC} \cdot \frac{\sin \angle BAM}{\sin \angle CAM}$$

and so

$$\frac{\sin \angle BAM}{\sin \angle CAM} = \frac{AC}{AB}$$

To sum up we obtain

$$\frac{\sin \angle CAX}{\sin \angle BAX} = \frac{\sin BAM}{\sin CAM}$$

and so the second condition of our Lemma is fulfilled, and we are done ∎

Here is our first example, which demonstrates the *thought process* on applying the Lemma.

> **Example 8: Russia 2011**
>
> Let $N$ be the midpoint of arc $ABC$ of the circumcircle of triangle $ABC$ and let $M$ be the midpoint of $AC$. Suppose that $I_1, I_2$ are the incenters of triangles $ABM$ and $CBM$. Prove that points $I_1, I_2, B, N$ lie on a circle.

*Proof:* It suffices to prove that $\angle I_1BI_2 = \angle I_1NI_2$, or equivalently that $\angle I_1NI_2 = \angle ANM$, which is equivalent to either proving that $\angle I_1NA = \angle MNI_2$, or proving that $\angle I_1NM = \angle I_2NC$. The key is to note that

$$\angle I_1NA + \angle I_1NM = \angle ANM = \angle MNC = \angle MNI_2 + \angle I_2NC,$$

and since $\angle ANC = \angle ABC < \pi$, by Lemma 4 it suffices to prove that

$$\frac{\sin \angle I_1NA}{\sin \angle I_1NM} = \frac{\sin \angle MNI_2}{\sin \angle I_2NC}$$

Now, note that

$$\angle I_1AM = \angle A/2 = 90° - \angle B/2 - \angle C/2 = \angle NCA - \angle I_2CA = \angle I_2CN$$

and similarly
$$\angle I_1AN = \angle I_2CM.$$

Moreover,
$$\angle I_2MN = \angle NMC - \angle I_2MC = 90° - \angle BMC/2 = \angle BMA/2 = \angle AMI_1$$

and similarly
$$\angle NMI_1 = \angle CMI_2.$$

Now, in order to finish, we may simply apply the Trigonometric version of Ceva's Theorem in triangles $ANM$ and $MNC$ (see the section "The Main weapons, Part 3") ∎

Our next example comes from the Mediterranean Olympiad.

> **Example 9: Mediterranean MO 1998**
>
> In a triangle $ABC$, $I$ is the incenter and $D, E, F$ are the points of tangency of the incircle with $BC, CA, AB$, respectively. The bisector of angle $\angle BIC$ meets $BC$ at point $M$, and the bisector of angle $\angle FDE$ meets $FE$ at point $P$. Prove that points $A, P, M$ are collinear.

*Proof:* We need to either prove that $\angle PAF = \angle MAB$ or that $\angle PAE = \angle MAC$. Since
$$\angle PAF + \angle PAE = \angle A = \angle MAB + \angle MAC,$$

from Lemma 4 it suffices to prove that
$$\frac{\sin \angle PAE}{\sin \angle PAF} = \frac{\sin \angle MAC}{\sin \angle MAB}$$

The key is to constantly apply Lemma 1. Note that
$$\frac{DE}{DF} = \frac{PE}{PF} = \frac{AE}{AF} \cdot \frac{\sin \angle PAE}{\sin \angle PAF},$$

and so
$$\frac{\sin \angle PAE}{\sin \angle PAF} = \frac{DE}{DF}$$

Moreover, on a similar manner we infer that
$$\frac{\sin \angle MAC}{\sin \angle MAB} = \frac{IC}{IB} \cdot \frac{AB}{AC},$$

and so we are left to prove that
$$\frac{IC}{IB} \cdot \frac{AB}{AC} = \frac{DE}{DF},$$

which boils down to a few applications of Theorem 1 to triangles $BFD, CDE$ and $BIC$, which are left to the reader to execute ∎

We end this section with an example from an old Mathlinks Contest.

> **Example 10: Mathlinks Contest 2008**
>
> Let $\Omega$ be the circumcircle of triangle $ABC$. Let $D$ be the point at which the incircle of $ABC$ touches its side $BC$. Let $M$ be the point on $\Omega$ such that the line $AM$ is parallel to $BC$. Also, let $P$ be the point at which the circle tangent to the segments $AB$ and $AC$ and to the circle $\Omega$ touches $\Omega$. Prove that the points $P, D, M$ are collinear.

*Proof (due to Silouanos Brazitikos):* We freely use the properties from the first Practice Problem for this Section. From the angle bisector theorem we have $\dfrac{QB}{QA} = \dfrac{PA}{PB}$ and $\dfrac{RA}{RC} = \dfrac{PA}{PC}$, and so dividing these two we have that

$$\frac{PC}{PB} = \frac{RC}{QB} \quad (1)$$

From the Law of Sines at the triangles $PBD, PDC$ we find that

$$\frac{\sin \angle BPD}{\sin \angle DPC} = \frac{CP \cdot BD}{BP \cdot CD} \quad (2)$$

Now we have that $\angle ARI = 90 - \dfrac{\angle A}{2}$ and $\angle RCI = \dfrac{\angle C}{2}$, hence $\angle RIC = \dfrac{\angle B}{2}$, so from the Law of Sines at the triangle $RIC$ we have that

$$RC = IC \cdot \frac{\sin \dfrac{\angle B}{2}}{\cos \dfrac{\angle A}{2}}$$

and similarly we have

$$BQ = IB \cdot \frac{\sin \dfrac{\angle C}{2}}{\cos \dfrac{\angle A}{2}}$$

Thus, $\dfrac{RC}{QB} = \dfrac{IC}{IB} \cdot \dfrac{\sin \dfrac{\angle B}{2}}{\sin \dfrac{\angle C}{2}}$, which combined with $\dfrac{IC}{IB} = \dfrac{CD}{BD} \cdot \dfrac{\cos \dfrac{\angle B}{2}}{\cos \dfrac{\angle C}{2}}$ implies that

$$\frac{RC}{QB} = \frac{CD}{BD} \cdot \frac{\sin \angle B}{\sin \angle C} \quad (3)$$

Combining relation (2) with relations (1) and (3), we may conclude that
$$\frac{\sin \angle BPD}{\sin \angle DPC} = \frac{\sin \angle B}{\sin \angle C}$$
However, by Lemma 4 this implies that $\angle DPC = \angle C$ and $\angle BPD = \angle B$. Now if we let line $PD$ meet the circumcircle at point $M'$, we will have that $\angle M'PC = C$, that is $\angle M'AC = C$, so the lines $AM'$ and $BC$ are parallel and points $M$ and $M'$ coincide, as desired ∎

# 5 Practice Problems

> **Problem 13: Miscellaneous (*)**
>
> Let $ABC$ be a triangle, and let $(\Gamma)$ be the circle tangent to $AB, AC$ and internally at $(ABC)$. Let $K, L, M$ be the tangency points of $(\Gamma)$ with $AB, AC$ and $(\Gamma)$, respectively, and $I$ be the incenter of $ABC$. Prove, using trigonometry, that
>
> (i) $I$ is the midpoint of $KL$ and
>
> (ii) Line $KI$ passes through the midpoint of arc $(BAC)$ in circle $(ABC)$.

> **Problem 14: Adapted from the Internet (**)**
>
> Let $ABC$ be a triangle and let its altitudes $AD, BE, CF$ concur at point $H$. Let $M$ be the midpoint of $EF$ and let lines $FE, BC$ intersect at point $P$. Assume that $PH$ meets $AB, AC$ at points $K, L$. Then, prove that points $A, M, Q$ are collinear, with point $Q$ being the intersection of the perpendiculars from points $K, L$ to $AB, AC$, respectively.

> **Problem 15: MEMO 2019/T5 (**)**
>
> Let $ABC$ be an acute-angled triangle such that $AB < AC$. Let $D$ be the point of intersection of the perpendicular bisector of the side $BC$ with the side $AC$. Let $P$ be a point on the shorter arc $AC$ of the circumcircle of the triangle $ABC$ such that $DP \parallel BC$. Finally, let $M$ be the midpoint of the side $AB$. Prove that $\angle APD = \angle MPB$.

> **Problem 16: Poland 2019/1 (**)**
>
> Let $ABC$ be an acute triangle. Points $X$ and $Y$ lie on the segments $AB$ and $AC$, respectively, such that $AX = AY$ and the segment $XY$ passes through the orthocenter of the triangle $ABC$. Lines tangent to the circumcircle of the triangle $AXY$ at points $X$ and $Y$ intersect at point $P$. Prove that points $A, B, C, P$ are concyclic.

> **Problem 17: Saint Petersburg 2022**
>
> In an acute triangle $ABC$ let $AH$ be the altitude from $A$ and $D$ be the antipode of $A$ on $(ABC)$. Point $I$ is the incenter of this triangle. Prove that $\angle BIH = \angle DIC$.

**Problem 18: Adapted from the Internet (**)**

Let $ABC$ be a triangle and let $I$ be its incenter. Let $D$ be the point of intersection of $AI$ with $BC$. Let $E$ be the incenter of triangle $ABD$ and let $F$ be the incenter of triangle $ADC$. Let line $DE$ intersect the circumcircle of triangle $BCE$ at point P ($\neq E$) and let line $DF$ intersect the circumcircle of triangle $BCF$ at point $Q$ ($\neq F$). Prove that the midpoint of $BC$ belongs to the circumcircle of triangle $DPQ$.

**Problem 19: Asian Pacific MO 2021 (***)**

Let $ABCD$ be a cyclic convex quadrilateral and $\Gamma$ be its circumcircle. Let $E$ be the intersection of the diagonals of $AC$ and $BD$. Let $L$ be the center of the circle tangent to sides $AB$, $BC$, and $CD$, and let $M$ be the midpoint of the arc $BC$ of $\Gamma$ not containing $A$ and $D$. Prove that the excenter of triangle $BCE$ opposite $E$ lies on the line $LM$.

**Problem 20: DGO 2022 (Individual) (***)**

Let $\omega$ denote the circumcircle of a scalene triangle $ABC$ with $AB < AC$. Tangents to $\omega$ at $B, C$ meet at $T$. Let $AT$ intersect $\omega$ again at $K$. The circumcircle of triangles $ACT, BKT$ at $T$ and $X$. The circumcircle of triangles $ABT, CKT$ at $T$ and $Y$. $BY$ intersects $CX$ at $Z$. Show that $AZ$ intersects the circumcircle of triangle $BTC$ at the midpoint of segment $XY$.

## 6 The main weapons, Part 3

In this section, we will introduce Menelaus's theorem, Ceva's theorem and - most importantly - it's trigonometric version.

> **Theorem 4: Menelaus's Theorem**
>
> Let $ABC$ be a triangle and $D, E, F$ points on $BC, CA, AB$, with $D$ on the extension of segment $BC$ and points $E, F$ in the interior of $AC, AB$, respectively. Then,
> $$D, E, F \text{ collinear} \iff \frac{DB}{DC} \cdot \frac{EC}{EA} \cdot \frac{FA}{FB} = 1$$

*Proof:* We only prove the right direction, leaving the converse as an exercise for the reader. If $X, Y, Z$ are on line $\overline{DEF}$ such that $AX, BY$ and $CZ$ are perpendicular to $\overline{DEF}$, then

$$\frac{DB}{DC} \cdot \frac{EC}{EA} \cdot \frac{FA}{FB} = \frac{BY}{CZ} \cdot \frac{CZ}{AX} \cdot \frac{AX}{BY} = 1,$$

as desired ∎

**Exercise 5:** Use the Ratio Lemma (Lemma 1) to provide one more proof of Menelaus's theorem.

> **Theorem 5: Ceva's theorem**
>
> Let $ABC$ be a triangle and $D, E, F$ points on $BC, CA, AB$ respectively. Then,
> $$AD, BE, CF \text{ concurrent} \iff \frac{DB}{DC} \cdot \frac{EC}{EA} \cdot \frac{FA}{FB} = 1$$

*Proof:* We present a proof based on areas (a section on which will be devoted on Part II of this paper). Let $[X]$ denote the are of shape $X$. Suppose that $AD, BE, CE$ concur at point $K$. Note that

$$\frac{DB}{DC} = \frac{[ADB]}{[ADC]} = \frac{[ODB]}{[ODC]} = \frac{[ADB] - [ODB]}{[ADC] - [ODC]} = \frac{[ABO]}{[ACO]}$$

Obtaining similar relations for the other ratios and multiplying, we conclude ∎

Here is a useful corollary of Ceva's theorem:

> **Lemma 6: One of the cevians is a median**
>
> Let $ABC$ be a triangle and $M$ be the midpoint of $BC$. If $E, F$ are points on $AB, AC$, respectively, then:
>
> $$EF \parallel BC \iff AM, BF, CE \text{ are collinear}$$

*Proof:* This is a direct application of Ceva's and Thales's theorem ∎

Here are two examples illustrating the importance of this Lemma.

> **Example 11: Balkan MO 2017/2**
>
> Consider an acute-angled triangle $ABC$ with $AB < AC$ and let $\omega$ be its circumscribed circle. Let $t_B$ and $t_C$ be the tangents to the circle $\omega$ at points $B$ and $C$, respectively, and let $L$ be their intersection. The straight line passing through the point $B$ and parallel to $AC$ intersects $t_C$ in point $D$. The straight line passing through the point $C$ and parallel to $AB$ intersects $t_B$ in point $E$. The circumcircle of the triangle $BDC$ intersects $AC$ in $T$, where $T$ is located between $A$ and $C$. The circumcircle of the triangle $BEC$ intersects the line $AB$ (or its extension) in $S$, where $B$ is located between $S$ and $A$. Prove that $ST$, $AL$, and $BC$ are concurrent.

*Proof (due to Dionysios Adamopoulos):* We will prove that $AB, AC$ are tangent to the circumcircles of triangles $BCD$ and $BCE$, respectively.

Note that $\angle SBD = \angle BCD$, and so it suffices to prove that

$$\angle DBA = 180° - \angle BCL \iff \angle DBC + \angle CBA =$$
$$= 180° - \angle BAC \iff \angle BCA + \angle CBA = 180° - \angle BAC,$$

which is true.

We work similarly for the second tangency. Moreover, $BT$ is antiparallel with $BC$ in triangle $ABC$. Hence the symmedian of vertex $A$, which is line $AL$ passes through the midpoint of $M$ of $BT$. Similarly $AL$ passes through the midpoint $N$ of $CS$.

However, note that $BT$ and $SC$ are parallel, due to an easy angle-chase. From the above Lemma in triangle $ASC$, we infer that $ST, BC$ and $AN \equiv AL$ concur, as desired ∎

> **Example 12: Adapted from the Internet**
>
> From point $S$, we draw the tangents $SA, SB$ to circle $(O, R)$. Let $M$ be a point on segment $SB$ and $N$ be the midpoint of $AM$. If $SN$ intersects $AB$ at point $P$, prove that $OP$ is perpendicular to $AM$.

*Proof:* Let $C$ be the point of intersection of $AM$ with $SO$ and $K$ be the intersection of $AB$ with $SO$. Points $S, O$ are equidistant from $A, B$, and so $OS$ is perpendicular to $AB$, hence $AK$ is perpendicular to $CO$..

In triangle $BAM$, points $K, N$ are the midpoints of $AB, AM$ respectively, hence:
$$\angle AKN = \angle ABM = \angle BAM$$

Now assume that $KN$ intersects $AS$ at point $T$. Note that
$$\angle AKT = \angle KAT, \ AT = TK = TS,$$

and so $T$ is the midpoint of $AS$.

From the above Lemma in triangle $KAS$, we obtain that $PC$ is parallel to $AS$, and since $AS$ is perpnedicular to $AO$, we conclude that $PC$ is perpendicular to $AO$.

From the above results, $P$ is the orthocenter of triangle $AOC$, finishing the problem ∎

Note that Ceva's theorem continues to hold when lines $AD, BE, CF$ concur at the *exterior* of triangle $ABC$. We will later see an example on this.

Here is the trigonometric version of Ceva's theorem.

> **Theorem 6: Trigononometric Ceva**
>
> Let $ABC$ be a triangle and $D, E, F$ points on sides $BC, CA, AB$ respectively. Then,
> $$AD, BE, CF \text{ concurrent} \iff \frac{\sin \angle DAB}{\sin \angle DAC} \cdot \frac{\sin \angle EBC}{\sin \angle EBA} \cdot \frac{\sin \angle FCA}{\sin \angle FCB} = 1$$

*Proof:* Apply Theorem 5 and Lemma 1 (Ratio Lemma) ∎

Three easy exercises follow, for the reader to try to implement how the above theorems work:

**Exercise 6:** Prove the existence of the orthocenter, the centroid, the incenter, the Nagel point, the Gergonne point and the Lemoine point using Ceva's theorem, as well as Trigonometric Ceva.

**Exercise 7:** Prove the existence of the isotomic and the isogonal conjugate of a point $P$ in the interior of a triangle $ABC$.

**Exercise 8 (due to Prodromos Fotiadis):** Consider a $n$–gon $A_1A_2...A_n$ and a point $K$ in its interior. Let $\theta_i = \angle KA_{i+1}A_i$ and $\varphi_i = KA_{i+1}A_{i+2}$ for all $i \in \{1, 2, \ldots, n\}$, with indices $\pmod{n}$ (i.e. $A_{n+1} \equiv A_1$ and $A_{n+2} \equiv A_2$). Prove that

$$\frac{\sin \angle \theta_1}{\sin \angle \varphi_1} \cdot \frac{\sin \angle \theta_2}{\sin \angle \varphi_2} \cdot \ldots \cdot \frac{\sin \angle \theta_n}{\sin \angle \varphi_n} = 1$$

Here is one more example.

> **Example 13: Adapted from the Internet**
>
> In an acute triangle $ABC$, such that the angle bisector $AD$, the median $BE$ and the altitude $CZ$ concur, prove that $\angle A \geq 45°$.

*Proof:* From Ceva's theorem,

$$\frac{AZ}{ZB}\frac{DB}{DC}\frac{EC}{AE} = 1$$

Using the angle bisector theorem and Lemma 6, we obtain $\frac{c}{b} = \frac{AE}{EC}$. However. $AE = c \cos A$ and $EC = a \cos C$, and so using the Law of Sines the previous relation rewrites as $\sin A \cos C = \sin B \cos A$, that is

$$\tan A = \frac{\sin B}{\cos C}$$

Since the triangle is acute, $\angle B + \angle C > \frac{\pi}{2}$, hence

$$\sin B > \sin(\pi/2 - \angle C) = \cos \angle C,$$

implying that $\tan A > 1$, which rewrites to what we wanted to prove ∎

In the spirit of the above easy example, we suggest the reader tries the following exercise, too:

**Exercise 9:** If in a triangle $ABC$ the median $AM$ the angle bisector $BD$ and the altitude $CE$ are equal and concur, can we infer that the triangle is *necessarily* equilateral?

We now present a useful theorem, and one application of it.

> **Theorem 7: Desargues's Theorem**
>
> Two triangles are in perspective axially if and only if they are in perspective centrally.

*Proof:* While projective solutions are much easier, as Desargues's Theorem is central in projective Geometry, we link to this solution using only Menelaus's theorem, from cut-the-knot ∎ Armed with this theorem, we present the following example.

> **Example 14: Adapted from the Internet**
>
> With diameters the sides $BC, CA, AB$ of right triangle $ABC$ (with $\angle A = 90°$), we draw the respective semicircles in its exterior, and let $M, N, L$ be their respective midpoints. Let now $S, P, Q$ be the intersections of lines $BA, MN$, $CA, ML$ and $CB, NL$, respectively. Prove that points $S, P, Q$ are collinear.

*Proof (due to Minos Margaritis):* We draw $CL$ and let point $E$ be the intersection of $CL, AB$ and point $E'$ be the midpoint of $AB$. Due to the similar triangles $ACE$ and $EE'C$ we obtain that

$$\frac{AE}{EE'} = \frac{AC}{E'L} = \frac{2AC}{AB} = 2\tan(\angle B)$$

Thus, we find that

$$\frac{AE}{EB} = \frac{\tan \angle B}{\tan \angle B + 1},$$

and similarly defining points $Z, Z' \in AC$, we obtain that

$$\frac{CZ}{ZA} = \frac{\tan \angle C + 1}{\tan \angle C}$$

Defining points $H, H' \in BC$, and since $AM$ is the angle bisector of $\angle A$ (as quadrilateral $ACMB$ is cyclic), we infer that

$$\frac{BH}{CH} = \frac{AB}{AC}$$

Note that

$$\frac{AE}{EB} \cdot \frac{CZ}{ZA} \cdot \frac{BH}{HC} = \frac{\tan \angle B}{\tan \angle B + 1} \cdot \frac{\tan \angle C + 1}{\tan \angle C} \cdot \frac{AB}{AC} = 1$$

(the proof of the last equality is an easy exercise), and so from the converse of Ceva's theorem in triangle $ABC$ for $CE, BZ, AH$ we obtain that these lines concur. This means that triangles $NML, ABC$ are perspective, and so from Desargues's theorem we obtain that points $S, P, Q$ are collinear, as desired ∎

Here is an easy exercise, which illustrates the importance of Desargues's theorem:

**Exercise 10:** Let $ABC$ be a triangle and let points $D, E$ belong on $AC, AB$, respectively. Let $BD, CE$ meet at point $O$ and let $M$ belong on line $OA$. Prove that lines $PQ, BC, DE$ concur, with $P$ being the intersection of $MD$ and $CE$, and $Q$ being the intersection of $ME$ and $BD$.

We finish this section with one more application of Menelaus's theorem, coupled with some elements from trigonometry.

> **Example 15: Adapted from the Internet**
>
> Let $ABC$ be a triangle inscribed in circle $(O)$ and let $M, N$ be the midpoints of $AB, AC$ and $Q, R$ be the feet of the altitudes from $B$ and $C$, respectively. Prove that lines $MQ, RN$ intersect on the Euler Line of triangle $ABC$ (i.e. line $HO$, with $H$ the orthocenter and $O$ the circumcenter).

*Proof:* Let $RN$ intersect $HO$ at point $S$ and $MQ$ intersect $HO$ at point $S'$. We ought to prove that points $S$ and $S'$ coincide. Moreover, let $HO$ intersect $AB, AC$ at points $X, Y$, respectively. From Menelaus's theorem in triangle $AXY$ we obtain that

$$\frac{S'X}{S'Y} = \frac{RX}{RA} \cdot \frac{NA}{NY},$$

and from Menelaus's theorem in triangle $AXY$ again we obtain that

$$\frac{SX}{SY} = \frac{MX}{MA} \cdot \frac{QA}{QY}$$

Hence, it is enough to prove that

$$\frac{RX}{RA} \cdot \frac{NA}{NY} = \frac{MX}{MA} \cdot \frac{QA}{QY}$$

Equivalently, the above relation rewrites as

$$\frac{RA}{MA} \cdot \frac{MX}{RX} = \frac{NA}{QA} \cdot \frac{QY}{NY}$$

Note that $RH$ and $OM$ are perpendicular to $AB$, and so are parallel. Thus,

$$\frac{MX}{RX} = \frac{MO}{RH}$$

and similarly

$$\frac{QY}{NY} = \frac{QH}{NO}$$

Putting it all together, it suffices to prove that

$$\frac{RA}{MA} \cdot \frac{MO}{RH} = \frac{NA}{QA} \cdot \frac{QH}{NO}$$

However,

- $\dfrac{RA}{MA} \cdot \dfrac{MO}{RH} = \dfrac{RA}{RH} \cdot \dfrac{MO}{MA} = \tan \angle AHR \cdot \tan \angle MAO = \tan \angle B \cdot \cot \angle C$ and

- $\dfrac{NA}{QA} \cdot \dfrac{QH}{NO} = \dfrac{NA}{NO} \cdot \dfrac{QH}{QA} = \tan \angle AON \cdot \tan \angle QAH = \tan \angle B \cdot \cot \angle C$

From the above bullets, we may conclude ∎

# 7 Practice Problems

**Problem 21: CJMO 2021 (*)**

Let $ABC$ be an acute triangle, and let the feet of the altitudes from $A$, $B$, $C$ to $\overline{BC}$, $\overline{CA}$, $\overline{AB}$ be $D$, $E$, $F$, respectively. Points $X \neq F$ and $Y \neq E$ lie on lines $CF$ and $BE$ respectively such that $\angle XAD = \angle DAB$ and $\angle YAD = \angle DAC$. Prove that $X, D, Y$ are collinear.

**Problem 22: Adapted from the Internet (*)**

Let triangle $ABC$ with $AB < AC$ have circumrcircle $(\Gamma)$ and circumcenter $O$. The tangents to $(\Gamma)$ to points $B$ and $C$ intersect at point $K$ and $OK$ intersects $BC$ at $Z$. Line $BO$ intersects lines $AZ$ and $AC$ at points $M$ and $X$ respectively and line $MC$ intersects $AB$ at point $Y$. Prove that the circumcircle of triangle $AXY$ and circle $(\Gamma)$ are tangent.

**Problem 23: Korea 1997 (*)**

In an acute triangle $ABC$ such that $AB \neq AC$, let $V$ be the intersection of the angle bisector of $A$ with $BC$, and let $D$ be the foot of the perpendicular from $A$ to $BC$. If $E$ and $F$ are the intersections of the cirucmcircle of triangle $AVD$ with $CA$ and $AB$ respectively, then show that lines $AD, BE, CF$ concur.

**Problem 24: Adapted from the Internet (*)**

Let $ABC$ be a quadrilateral inscribed at circle $(O)$. Diagonals $AC, BD$ intersect at point $K$. If lines $AB, DC$ intersect at point $M$ and the tangents to $(O)$ at points $B$ and $C$ at point $L$, then prove that points $K, L, M$ are collinear.

**Problem 25: IMO Shortlist 2001 (**)**

Let $A_1$ be the center of the square inscribed in acute triangle $ABC$ with two vertices of the square on side $BC$. Thus one of the two remaining vertices of the square is on side $AB$ and the other is on $AC$. Points $B_1$, $C_1$ are defined in a similar way for inscribed squares with two vertices on sides $AC$ and $AB$, respectively. Prove that lines $AA_1$, $BB_1$, $CC_1$ are concurrent.

**Problem 26: Adapted from the Internet (**)**

Let triangle $ABC$ be such that $AB + AC = 2BC$. Its incircle touches $AB, AC$ at points $D, E$, respectively. Prove that line $DE$ bisects the median $AM$.

**Problem 27: USAMO 2008**

Let $ABC$ be an acute, scalene triangle, and let $M$, $N$, and $P$ be the midpoints of $\overline{BC}, \overline{CA}$, and $\overline{AB}$, respectively. Let the perpendicular bisectors of $\overline{AB}$ and $\overline{AC}$ intersect ray $AM$ in points $D$ and $E$ respectively, and let lines $BD$ and $CE$ intersect in point $F$, inside of triangle $ABC$. Prove that points $A$, $N$, $F$, and $P$ all lie on one circle.

**Problem 28: Adapted from the Internet (**)**

Triangle $ABC$ has $AB = AC$. Let point $E$ be in its interior, belonging at the perpendicular bisector of $AB$. On line $AE$ let point $Z$ be such that $\angle ABZ = \angle BCE$.
If $K$ is the circumcenter of triangle $ABZ$ and $H$ the reflection of $K$ across $BZ$, then prove that $HB = HC$.

**Problem 29: Romania TST 2002 (**)**

Let $M$ and $N$ be the midpoints of the respective sides $AB$ and $AC$ of an acute-angled triangle $ABC$. Let $P$ be the foot of the perpendicular from $N$ onto $BC$ and let $A_1$ be the midpoint of $MP$. Points $B_1$ and $C_1$ are obtained similarly. If $AA_1$, $BB_1$ and $CC_1$ are concurrent, show that the triangle $ABC$ is isosceles.

**Problem 30: The Cevian Nest Theorem (**)**

Let $D, E, F$ be points on sides $BC, CA, AB$ respectively of a triangle $ABC$. Also, let $X, Y, Z$ be points on the sides $EF, FD, DE$ respectively of triangle $DEF$. Consider the three triples of lines $(AX, BY, CZ), (AD, BE, CF)$ and $(DX, EY, FZ)$. Prove that if any two of these triples are concurrent, then the thirs one is as well.

**Problem 31: International Olympiad of Metropolises 2020 (**)**

In convex pentagon $ABCDE$ points $A_1, B_1, C_1, D_1, E_1$ are intersections of pairs of diagonals $(BD, CE)$, $(CE, DA)$, $(DA, EB)$, $(EB, AC)$ and $(AC, BD)$ respectively. Prove that if four of quadrilaterals $AB_1A_1B$, $BC_1B_1C$, $CD_1C_1D$, $DE_1D_1E$ and $EA_1E_1A$ are cyclic then the fifth one is also cyclic.

(*Remark:* Although no solution using Menelaus's or Ceva's theorem is known, the ideas used in them can provide sufficient motivation for solving this problem.)

**Problem 32: Ukraine IMO TST (***)**

Let altitudes $AH_1$ and $BH_2$ of an acute triangle $ABC$ intersect at point $H$. Let $\gamma_1$ be the circle passing through $H_2$ and tangent to $BC$ at $H_1$ and let $\gamma_2$ be the circle passing through $H_1$ and tangent to $AC$ at point $H_2$. Prove that the second tangent $BX$ of circle $\gamma_1$ and the second tangent $AY$ of circle $\gamma_2$ (with $X \neq H_1$ and $Y \neq H_2$) intersect on the circumrcircle of triangle $XHY$.

# 8 Historical information

In this section, we will provide short biographies of all the mathematicians whose names were involved in Theorems or Lemmas in this paper. The mathematicians are listed in chronological order. Shedding light into the background of the mathematical explorations often provides us with useful insight, and helps to comprehend the perpetual nature of Mathematics.

1. **Menelaus of Alexandria (c. 70 - 140 AD):** A prominent Greek Mathematician and Astronomer, he was raised in Alexandria and it is suspected that he later moved to Rome, where he spent the remaining of his life. He is credited with introducing the concepts of *spherical triangle* and *geodesics*. His most important known work is *Sphaerica*, which is composed of three books and provides useful applications of the geometry of the sphere on astronomical calculations.

2. **Girard Desargues (1591 - 1661):** A French mathematician and engineer, widely considered as one of the founders of *projective geometry*. Desargues's theorem is fundamendal in projective geometry, and even leads to the classification of planes in *Desarguesian* and *non - Desarguesian*, depending on whether the theorem applies on these planes, or not. Desargues worked as a private tutor, engineer and architect, and made contributions in engineering, such as the designation of a system for raising water that he installed near Paris. His work in projective geometry is deeply influenced by that of Johannes Kepler in optics, and other scientists of the medieval age. Desargues deceased in September 1661, in Lyon.

3. **Giovanni Ceva (1647 - 1734):** Italian Mathematician, known for proving the relevant theorem in elementary geometry. Born in Milan, he received his education there and later studied at the University of Pisa, in which he became a professor of Mathematics. He resigned from this position, though, and became the Professor of Mathematics at the University of Mantua, in which city he lived for the rest of his life. His work was mainly in elementary geometry. He is credited with proving Ceva's theorem (known since the 11th century), and rediscovering Menelaus's theorem. Ceva also worked in applied mathematics, studying *hydraulics* and *statics*. Ceva deceased in May 1734, in Mantua.

4. **Leonhard Euler (1707 - 1783):** Swiss mathematician, of the greatest of all time. His contributions are so huge and unique, that we are obliged to link to this document for relevant information regarding his life and his extraordinary work.

5. **Joseph Diez Gergonne (1771 - 1859):** French Mathematician and logician, Born in Nancy, he was a captain of the French army until 1795, with small breaks back to civilian life. From 1795 and on, he took the chair of "transcendental mathematics" at Nîmes. Following his speciality,

Geometry, Gergonne founded his own journal, *Annales de Mathématiques Pures et Appliquées* ("Annals of Pure and Applied Mathematics"), which he published for more than two decades (1810 - 1831). Many renowned mathematicians, such as Cauchy, Steiner, Brianchon, Poncelet, Liouville, Poisson, Ampere and Chasles contributed to Gergonne's magazine, which apart from geometry featured artcicles on history, philosophy, and mathematics education. Gergonne's work in Geometry consists of exploring the *principle of duality* in projective geometry and introducing the concept of a *polar*, while he also made contributions in Astronomy - he was appointed to the chair of astronomy at the University of Montpellier in 1816, at which city he also deceased in May 1859.

6. **Jakob Steiner (1796 - 1863):** Swiss mathematician, most known for his work on Geometry. Born in a small village in the Canton of Bern, he studied in Heidelberg and proceeded to earn his living by tutoring students in Berlin. His contributions in Geometry include introducing the *Steiner conic*, proving the *Steiner - Lehmus* theorem, introducing the notion of the *Steiner point*, stating the *Miquel and Steiner's quadrilateral theorem* (which was subsequently proved by Miquel), and discovering the *Roman (or Steiner) surface*. Steiner was a huge advocate of synthetic Geometry, strongly hating analytical tools. As a result, Steiner is widely known for his rigorous and beautiful proofs in the whole spectrum of geometry. Steiner also made some contributions in Combinatorics, introducing what is now known as *Steiner systems*. Steiner deceased in April 1863, in Bern.

7. **Christian Heinrich von Nagel (1803 - 1882):** Nagel was born in Stuttgart, and he initially studied theology at the Tübinger Stift. Soon he discovered his interest in Mathematics, and became the mathematics and science teacher at the local high school. He later moved to Ulm, where he also was a teacher at the local Gymnasium. His doctorate paper *De triangulis rectangulis ex algebraica aequatione construendis* ("About right triangles construable from an algebric equation") was the start of his geometry endeavors, which culminated in the naming of a very important triangle center by his name. Nagel decased in October 1882, in Ulm.

8. **Émile Lemoine (1840 - 1912):** French mathematician and civil engineer, whose full name was *Émile Michel Hyacinthe Lemoine*. Lemoine studied at a variety of schools, of which the most prominent is the École polytechnique - he later became a professor there for many years. For more than two decades he published the *L'Intermédiaire des mathématiciens* magazine, along with Charles-Ange Laisant. His main work interests were Geometry, civil engineering and - amatuerly - music. Lemoine's work in geometry has been so fundamental that he is credited as the co-founder (along with Brocard and Neuberg) of modern plane geometry. His most important works are proving the existence of the *Lemoine point and circle*, publishing the *Géométrographie*, a work which attempted to create a methodological system by which constructions could be judged, and writ-ing the *La Géométrographie ou l'art des constructions géométriques* treatise concerning compass and straightedge constructions. Lemoine also posed a relevant conjecture in Number Theory which still remains open up to this day and is considered of equal difficulty as Goldbach's conjecture. Lemoine deceased in February 1912, in Paris.

9. **Hiroshi Haruki (? - 1997):** Japanese mathematician, best known for his work in functional equations and his discoveries in plane geometry. He is credited with discovering *Haruki's theorem* and *Haruki's lemma* and with performing significant research on Conic Sections. Haruki worked at Osaka University in Japan and, for more than two decades, at the University of Waterloo, at which university he was a founding member of its computer science department. Haruki deceased in September 1997.

# 9 Abbreviations

Throughout the text, a number of abbreviations were used in citing the source of the problems proposed.

1. *Country, Year:* The problem is from a mathematical competition conducted in that specific country and that specific year, e.g. *Korea 1997*.
2. *USAMO:* United States of America Mathematical Olympiad.
3. *Balkan MO:* Balkan Mathematical Olympiad.
4. *TST:* Team Selection Test.
5. *TSTST:* Team Selection Test for the Selection Test.
6. *EGMO:* European Girls Mathematical Olympiad.
7. *OTJMO:* Online Tests Junior Mathematical Olympiad (student - run competition).
8. *CJMO:* Christmas Junior Mathematical Olympiad (student - run competition)
9. *GGG:* Competition run by Andrew Wu, see [here](#).
10. *Mediterranean MO:* Meditterranean Mathematics Contest.
11. *MEMO:* Middle European Mathematical Olympiad.
12. *DGO:* Discord Geometry Olympiad (student - run competition).